\documentclass[11pt, a4paper]{amsart}

\usepackage{amsmath,amsthm,amssymb}
\usepackage[foot]{amsaddr}

\newtheorem{rema}{Remark}
\newcommand{\br}{\begin{rem}}
\newcommand{\er}{\end{rem}}
\newcommand{\bex}{\begin{exa}}
\newcommand{\eex}{\end{exa}}
\newcommand{\bd}{\begin{Def}}
\newcommand{\ed}{\end{Def}}
\newtheorem{theorem}[rema]{Theorem}
\newcommand{\bt}{\begin{theorem}}
\newcommand{\et}{\end{theorem}}
\newtheorem{lemma}[rema]{Lemma}

\newcommand{\bl}{\begin{lemma}}
\newcommand{\el}{\end{lemma}}
\newcommand{\be}{\begin{equation}}
\newcommand{\ee}{\end{equation}}
\newcommand{\bea}{\begin{eqnarray}}
\newcommand{\eea}{\end{eqnarray}}

\newcommand{\adots}{\mathinner{\mkern2mu\raise1pt\hbox{.}\mkern2mu
\raise4pt\hbox{.}\mkern2mu\raise7pt\hbox{.}\mkern1mu}}

\newcommand{\beq}{\begin{equation}}  
\newcommand{\eeq}{\end{equation}}  
\newcommand{\bear}{\begin{array}}  
\newcommand{\eear}{\end{array}} 
\newcommand\la{{\lambda}}   
\newcommand\La{{\Lambda}}  
\newcommand\De{{\Delta}}   
   
\newcommand\al{{\alpha}}

\newcommand\eps{{\epsilon}}  
    
\newcommand\mma{{\bf M}}

\newtheorem{thm}{Theorem}[section]

\newtheorem{rem}[thm]{Remark}

\newtheorem{exa}[thm]{Example}

\newtheorem{lem}[thm]{Lemma}

\newtheorem{cor}[thm]{Corollary} 

\newenvironment{prf}{\trivlist \item [\hskip 
\labelsep {\bf Proof:}]\ignorespaces}{\qed \endtrivlist} 
 
\theoremstyle{remark}

\newcommand{\Z}{\mathbb{Z}}
\newcommand{\Zo}{\mathbb{Z}_{\ge0}}
\newcommand{\Zp}{\mathbb{Z}_{>0}}
\newcommand{\Q}{\mathbb{Q}}

\newcommand{\R}{{\mathbb R}}

\title[Continued fractions and irrationality exponents]{Continued fractions and irrationality exponents 
for modified Engel and Pierce series}
\author{Andrew N. W. Hone}
\address[A. N. W. Hone]{School of Mathematics, Statistics and Actuarial Science, University of Kent,
Canterbury CT2 7FS, UK. Currently on leave at School of Mathematics and Statistics, 
University of New South Wales, Sydney, NSW 2052, Australia.}
\email[A. N. W. Hone]{anwh@kent.ac.uk}

\author{Juan Luis Varona}
\address[J. L. Varona]{Departamento de Matem\'aticas y Computaci\'on, 
Universidad de La Rioja, 26006 Logro\~no, Spain}
\email[J. L. Varona]{jvarona@unirioja.es}



\begin{document}
%----------------

\begin{abstract} An Engel series is a 
 sum of %the 
reciprocals of a non-decreasing sequence $(x_n)$ of positive integers,
 which is such that each term is divisible by the previous one, and a Pierce series 
is an alternating sum of the reciprocals of a sequence with the same 
property. Given an arbitrary rational number, 
we show that 
there is a family of Engel series which when added to it  
produces a transcendental number $\al$ whose continued fraction expansion is 
determined explicitly by % in terms of $r$, $m$ and 
the corresponding 
sequence $(x_n)$, where the latter is generated 
by a certain  nonlinear recurrence of second order. 
We also present an analogous result for a rational number 
with a Pierce series added to or subtracted from it. In both situations (a 
rational number combined with either an Engel or a Pierce series), 
the irrationality exponent is bounded 
below by $(3+\sqrt{5})/2$, and we further identify 
infinite families of transcendental numbers $\al$ whose 
irrationality exponent can be computed precisely.   
We also construct the continued fraction expansion 
for an arbitrary rational number added to 
an Engel series with the stronger 
property that $x_j^2$ divides $x_{j+1}$ for all $j$.
\end{abstract}

%----------------
\maketitle
%----------------

%----------------
\section{Introduction} %and main results}
\label{sec:intro}
%----------------

Given a sequence of positive integers $(x_n)$, which is such that $x_n | x_{n+1}$ %$x_n \mid x_{n+1}$
 for all $n$, the sum of the reciprocals is the Engel series
\begin{equation}
\label{eq:Engel}
  \sum_{j=1}^{\infty} \frac{1}{x_j}
  = \sum_{j=1}^{\infty} \frac{1}{y_1y_2 \cdots y_j},
\end{equation}
where $y_1 = x_1$ and $y_{n+1} = x_{n+1}/x_n$ for $n \ge 1$,
and the alternating sum of the reciprocals is the Pierce series
\begin{equation}
\label{eq:Pierce}
  \sum_{j=1}^{\infty} \frac{(-1)^{j+1}}{x_j}
  = \sum_{j=1}^{\infty} \frac{(-1)^{j+1}}{y_1y_2 \cdots y_j}.
\end{equation}
(It should be assumed that $(x_n)$ is eventually increasing, in the 
sense that for all $n$ there is some $n'>n$ with $x_{n'}>x_n$, which guarantees the convergence of both sums~\eqref{eq:Engel} and~\eqref{eq:Pierce}.) 
Every positive real number admits both an Engel expansion, of the 
form \eqref{eq:Engel}, and a Pierce expansion \eqref{eq:Pierce} \cite{duverney}. Although they are not quite so well 
known, Engel expansions and Pierce expansions are in many ways analogous to 
continued fraction expansions, both in the sense that they are determined recursively, and from a 
metrical point of view; for instance, see \cite{erdosetal} for the case of Engel series. 

In recent work \cite{Hone1}, %, Hone2}, A.N.W. Hone 
the first author presented a family of sequences $(x_n)$ generated by a nonlinear recurrences of second order, 
of the form 
\beq\label{2nd} 
x_{n+1}x_{n-1} = x_n^2 \,\Big(1+x_n\,G(x_n)\Big), \qquad n\geq 2,  \qquad G(x)\in \Z [x], 
\eeq 
where the polynomial  $G$  takes positive values at positive arguments, such 
that the corresponding Engel series~\eqref{eq:Engel} yields a transcendental number 
whose 
continued fraction expansion is explicitly given in terms of  the $x_n$. More recently \cite{Var2}, the second 
author proved that, when the sequence $(x_n)$ is generated by a recurrence like 
(\ref{2nd}),  an analogous result holds for the associated Pierce series~\eqref{eq:Pierce}, 
although the structure of the corresponding continued fractions is different. In fact, in the latter 
work the polynomial $G(x_n)$  was replaced by an arbitrary sequence of 
positive integers,  as it had already been noted in \cite{Hone3} that the recurrence 
(\ref{2nd}) could be modified in this way and further allow the 
explicit continued fraction expansion to be determined for the sum 
of an arbitrary rational number $r=p/q$ and an Engel series, that is 
\beq\label{pqengel} 
\al = \frac{p}{q} + \sum_{j=2}^\infty  \frac{1}{x_j}, \qquad 
\mathrm{with}\,\,\, x_1=q.
\eeq   

In the next section we show that 
the initial conditions for the sequence $(x_n)$ can be specified 
more generally than  in %compared with 
\cite{Hone3},   
allowing dependence on a non-negative integer parameter $m$, 
and present analogous 
results for a  
family of transcendental numbers defined by a rational number with 
a Pierce series added to or subtracted from it, of the form 
\beq\label{pqpierce} 
\al = \frac{p}{q} \pm \sum_{j=2}^\infty  \frac{(-1)^j}{x_j}.
\eeq 
%where once again the sequence $(x_n)$ is generated 
%by a nonlinear recurrence of second order that generalizes (\ref{2nd}). 
%Furthermore, 
In section 3 it is proved that, 
in both cases  (\ref{pqengel}) and (\ref{pqpierce}),  %we show  that 
$\al$ has 
irrationality exponent $\mu(\al)\geq (3+\sqrt{5})/2$, and 
if the 
nonlinear recurrence for $(x_n)$  has a particular form then 
%the irrationality exponent of 
$\mu(\al)$ can be computed precisely. 
Explicit continued fractions for series of the form 
 (\ref{pqengel})  
with the stronger property that $x_j^2|x_{j+1}$
for all $j$ are constructed in 
the final section, generalizing the results in \cite{Hone2}. 

\section{Explicit continued fractions} 
\setcounter{equation}{0}

Before proceeding, we fix our notation for continued fractions and briefly mention some of their  
standard properties,  which can be found in many books %, for instance, 
\cite{Cas, duverney, Khin}. 
We denote a finite continued fraction by  
\beq\label{ficf}
  [a_0;a_1,a_2,\dots,a_n]
  = a_0 + \cfrac{1}{ a_1 + \cfrac{1}{ a_2 + \cfrac{1}{ \cdots +\cfrac{1}{a_n}} } }
  = \frac{p_n}{q_n},
\eeq
where  $a_0 \in \Z$, $a_j \in \Zp$  and $p_n/q_n$ is in lowest terms with $q_n>0$. Every $r\in\Q$ 
can be written as a finite continued fraction (\ref{ficf}), although 
this representation is not unique (see (\ref{props}) below). Each $\al\in\R\setminus\Q$ 
is given uniquely by an infinite continued fraction with convergents $p_n/q_n$ 
of the form (\ref{ficf}), that is (with $a_0=\left \lfloor{\al}\right \rfloor $) 
\beq\label{al}
\al = 
  [a_0;a_1,a_2,\dots] = \lim_{n\to\infty} [a_0;a_1,a_2,\dots,a_n] = \lim_{n\to\infty}\frac{p_n}{q_n}.
\eeq 
%and this limit always exists. 

The three-term recurrence relation satisfied by the  numerators and denominators of the convergents 
%can be conveniently 
is encoded in the matrix relation 
\beq\label{mat} 
\left(\bear{cc} p_{n+1} & p_n \\ 
q_{n+1} & q_n 
\eear\right) = 
\left(\bear{cc} p_n &  p_{n-1}\\ 
q_n &  q_{n-1}
\eear\right)\left(\bear{cc} a_{n+1} & 1 \\ 
1 & 0
\eear\right), 
\eeq 
valid for $n\geq -1$, with
$$  
\left(\bear{cc} p_{-1} & p_{-2} \\ 
q_{-1} & q_{-2} 
\eear\right) = \left(\bear{cc} 1 & 0 \\ 
0 & 1
\eear\right)
.
$$
By taking the determinant of both sides of (\ref{mat}), one obtains the
identity  
\beq\label{detid} 
  p_j q_{j-1} - p_{j-1} q_j = (-1)^{j-1}, \qquad j \ge 1. 
\eeq 

Note that  any finite continued fraction can be rewritten as another one of 
different length, since   one can always apply one of two 
operations, namely 
\beq\label{props} 
  [a_0;a_1,a_2,\dots,a_{j}] =
\begin{cases}
  [a_0;a_1,a_2,\dots,a_{j}-1,1], & \text{if $j=0$ or $a_j >1$},\\
  [a_0;a_1,a_2,\dots,a_{j-1}+1], & \text{if $a_j = 1$},
\end{cases}
\eeq 
and both operations change the parity of the length.

Henceforth we  fix a rational number $r=p/q$ in lowest terms, with $q \ge 1$, and an integer parameter $m\in\Zo$. 
Without loss of 
generality, because of (\ref{props}), we may specify the continued fraction of $r$ to be 
\beq\label{evencf} 
\frac{p}{q}=  [a_0;a_1,a_2,\dots,a_{2k}], 
\eeq 
with the index of the last partial quotient being even 
(and shifting $r$ by an integer changes $a_0$ but 
otherwise makes no difference). Given %a positive integer 
$y_0\in\Zp$, which 
in due course will be fixed differently  according to the context, we define two sequences 
$(y_n)_{n\geq 0}$ and $(x_n)_{n\geq 1}$ via the recursion relations 
\beq\label{recrels} 
y_{j}=y_{j-1}(1+u_{j}x_{j}), \quad 
x_{j+1}= x_{j}y_j, \quad \mathrm{for} \,\,j\geq 1, \, \, \mathrm{with}\,\, x_1=q, 
\eeq  
where $(u_n)_{n\geq 1}$ is an arbitrary sequence of positive integers. Note that 
the second relation guarantees the property $x_j|x_{j+1}$ required 
for an Engel series or a Pierce series. It is an 
immediate consequence of (\ref{recrels}) that, given $x_1=q$ and 
$x_2=x_1y_1=qy_0(1+u_1q)$, the subsequent terms of the 
sequence $(x_n)$ are determined by the nonlinear recurrence 
\beq\label{nonl} 
x_{n+1}x_{n-1} = x_n^2 \,(1+u_n\, x_n), \qquad n\geq 2. 
\eeq  
It may happen that the sequence $(u_n)$ is defined entirely in terms of 
the sequence $(x_n)$, by specifying a function $G$ such that 
\beq\label{aut} 
u_n = G(x_n), \qquad G: \, \Zp\rightarrow\Zp, 
\eeq 
which means that (\ref{nonl}) becomes an autonomous recurrence of the 
form (\ref{2nd}), although the function $G$ need not necessarily be a polynomial. 
Yet in general the recurrence (\ref{nonl}) is non-autonomous, whenever % the sequence 
$(u_n)$ is specified by 
\beq\label{naut} 
u_n = \widehat{G}(x_n,n), \qquad \widehat{G}: \, \Zp\times\Zp\rightarrow\Zp, 
\eeq 
where $\widehat{G}$ is a non-trivial function of its second argument. 

The following theorem extends the results of \cite{Hone3} (corresponding to the case $m=0$), 
being themselves extensions of \cite{Hone1} (corresponding to $p/q=1$). 

\begin{thm} \label{engelthm}
Let $\al$ be given by the series  
\beq \label{aleng} 
\al = \frac{p}{q}+\sum_{j=2}^\infty \frac{1}{x_j}, 
\eeq 
for $(x_n)$ defined by (\ref{recrels}) with 
%generated from  the reciprocals of the terms of the  sequence $(x_n)$ defined by (\ref{recrels}) with 
\beq\label{y0valeng}
y_0=mq+q_{2k-1}+1, 
\eeq 
where $q_{2k-1}$ is the denominator of the $(2k-1)$th convergent of 
(\ref{evencf}). 
Then  the continued fraction expansion of $\al$ has the form (\ref{al}), 
where the partial quotients $a_j$ coincide with those of 
(\ref{evencf}) for $0\leq j\leq 2k$, while  
\beq\label{engcf}
a_{2k+1}=m+u_1 y_0, \quad and \quad
a_{2k+2j}=x_j, \, \, a_{2k+2j+1}=u_{j+1}y_j, \,\, %for \quad 
\forall j\geq 1.
\eeq 
\end{thm} 
\begin{prf} 
The proof  consists of showing that the partial sums of the series (\ref{aleng}) 
coincide with the convergents with even index, that is 
\beq\label{parts} 
\frac{p}{q}+\sum_{j=2}^{N}\frac{1}{x_j} = \frac{p_{2k+2N-2}}{q_{2k+2N-2}}, 
\eeq
and then taking the limit $N\to\infty$. We omit the details, 
since the inductive proof of 
(\ref{parts}) %by induction on $N$ 
is almost identical to that of Theorem 2.1 
in \cite{Hone3}, with part of the %inductive 
hypothesis being that the denominators 
of the convergents are %given by 
\beq\label{cvgtseng} 
q_{2k+2N-3}=y_{N-1}-1, \qquad q_{2k+2N-2}=x_N
\eeq  
for $N\geq 1$.  The only difference is in verifying the base step: 
%of the induction: 
specifically, %checking  
that, for $N=1$, the three-term recurrence 
$q_{2k+2N-1} = a_{2k+2N-1}q_{2k+2N-2} + q_{2k+2N-3}$ 
gives the correct expression for $q_{2k+1}$. But when $N=1$, 
by using (\ref{y0valeng}) and (\ref{engcf}), together with the first recursive 
relation in (\ref{recrels}), 
the right-hand side %of the above 
becomes 
$
(m+u_1 y_0)q+ y_0 - mq-1 = y_0(1+u_1x_1)-1=y_1-1$,
which is the required formula for $q_{2k+1}$. 
%All other steps of the proof are the same as in \cite{Hone3}.
\end{prf} 

Continued fractions for some alternating series whose
sum is a transcendental number were considered in \cite{DaSha}, 
and the Pierce series in  \cite{Var2} provide other examples. 
As a first attempt at generalizing the latter results, 
we consider a rational number added to a Pierce series.

\begin{thm}\label{pierce1} 
Suppose $x_1=q>1$, and let $\al$ be given by the series  
\beq \label{alpie1} 
\al = \frac{p}{q}+\sum_{j=2}^\infty \frac{(-1)^j}{x_j}, 
\eeq 
for $(x_n)$ defined by (\ref{recrels}) with 
\beq\label{y0valpie1}
y_0=mq+q_{2k-1}-1, 
\eeq 
where $q_{2k-1}$ is the denominator of the $(2k-1)$th convergent of 
(\ref{evencf}), and the restriction $m\geq 1$ should be imposed 
if $q_{2k-1}=1$. 
Then  the continued fraction expansion of $\al$ has the form (\ref{al}), 
where the partial quotients $a_j$ coincide with those of 
(\ref{evencf}) for $0\leq j\leq 2k$, while  
\beq\label{pie1cf} 
a_{2k+1}=m+u_1 y_0-1, 
\eeq 
and 
\beq \label{pie1coeff}
a_{2k+3j-1}=1, \, \, a_{2k+3j}=x_j-1, \, \, a_{2k+3j+1}=u_{j+1}y_j-1 \,\,\, %for \quad 
\forall j\geq 1.
\eeq 
\end{thm} 
\begin{prf} Let $S_N$ denote the $N$th partial sum of the series (\ref{alpie1}), that is
$$ 
S_N= \frac{p}{q}+\sum_{j=2}^N \frac{(-1)^j}{x_j}. 
$$
We will show by induction that 
$$ 
S_N=[a_0; a_1, a_2,\ldots, a_{2k+3N-3}],
$$
and then the result follows in the limit $N\to\infty$. 
As part of the inductive hypothesis, we also require 
the following expressions for the denominators 
of the convergents, for $N\geq 1$: 
\beq\label{cvgtspie1} 
q_{2k+3N-3}=x_{N}, \, q_{2k+3N-2}=y_N-x_N+1, 
\,   q_{2k+3N-1}=y_N+1 . 
\eeq   
For the base case $N=1$ it is clear that 
$S_1=p/q$ is given by  (\ref{evencf}), and 
$q_{2k}=q=x_1$, while by using the three-term recurrence 
it follows from (\ref{y0valpie1}), (\ref{pie1cf}) and 
(\ref{recrels}) that 
$$ \bear{rcl} q_{2k+1}=a_{2k+1}q_{2k}+ q_{2k-1}
& =& (m+u_1 y_0-1)q +y_0 -mq+1 \\ 
& = &  y_0(1+u_1x_1) -x_1+1  \\
& =& y_1-x_1+1,
\eear
$$ 
and 
$q_{2k+2}=a_{2k+2}q_{2k+1}+q_{2k} =y_1-x_1+1+x_1=y_1+1$, 
which confirms (\ref{cvgtspie1}) in this case. For the inductive step, 
the first expression in (\ref{cvgtspie1}) 
is also  verified with the three-term recurrence,  as 
$$ 
\bear{rcl}
q_{2k+3N}& = & a_{2k+3N}\, q_{2k+3N-1}+q_{2k+3N-2}\\
& =& (x_N-1)(y_N+1)+y_N-x_N+1\\ 
&=&x_Ny_N=x_{N+1}, 
\eear 
$$ 
using (\ref{recrels}) once more, 
and the other two expressions for the denominators of the 
convergents are verified similarly. 
Now let $\mma_{n+1}$ denote the matrix on the left-hand 
side of (\ref{mat}), and let 
$$ 
{\bf A}_n = \left( \bear{cc} a_n & 1 \\ 1 & 0 \eear\right); $$ 
then observe that
$$\bear{rcl}\mma_{2k+3N}& =& \mma_{2k+3N-3} \,
{\bf A}_{2k+3N-2}  \,
 {\bf A}_{2k+3N-1} \, {\bf A}_{2k+3N}  \\
& =& \mma_{2k+3N-3} 
\left(\bear{cc} u_N x_N y_{N-1}-1 & u_N y_{N-1} \\ x_N & 1 \eear\right) .
\eear 
$$ 
The first column of the above identity yields the formula 
\beq\label{pform} 
p_{2k+3N}=(u_N x_N y_{N-1}-1)p_{2k+3N-3} + x_N p_{2k+3N-4} ,
\eeq 
and the analogous expression for $q_{2k+3N}$. 
Hence 
$$ \bear{rcl}
p_{2k+3N}  
& = & (y_N- y_{N-1}-1)p_{2k+3N-3} + x_N p_{2k+3N-4} \\
& = & (y_N- q_{2k+3N-4})p_{2k+3N-3} + q_{3k+3N-3} p_{2k+3N-4} \\ 
&=& y_N \, p_{2k+3N-3} + (-1)^{N+1}, 
\eear
$$
where we have used (\ref{cvgtspie1}) and the first relation in (\ref{recrels}), 
followed by (\ref{detid}). (Note that the calculation leading to the latter 
expression is slightly different in the base case $N=1$, involving 
the use of  (\ref{y0valpie1}) and (\ref{pie1cf}), but the 
conclusion is the same.) Thus, using the first formula in (\ref{cvgtspie1}), 
and the fact that $y_N /x_{N+1}=1/x_N$, we have 
$$ \bear{rcl} {p_{2k+3N}}/{q_{2k+3N}} & = & 
{x_{N+1}}^{-1} \left(y_N \, p_{2k+3N-3} + (-1)^{N+1}\right) \\
&=& {p_{2k+3N-3}}/{q_{2k+3N-3}} + {(-1)^{N+1}}/{x_{N+1}}. \eear$$ 
So by the inductive hypothesis, 
$$ 
\frac{p_{2k+3N}}{q_{2k+3N}} = 
 S_N +  \frac{(-1)^{N+1}}{x_{N+1}} = S_{N+1}, 
$$ 
as required.
\end{prf}
\begin{rem} 
 The assumption that $q>1$ implies $k\geq 1$ in (\ref{evencf}), and is made 
to ensure that 
$a_{2k+3}>0$. However, when $q=1$, the appearance of a zero in the continued fraction 
can be dealt with by applying the concatenation operation 
\beq\label{concat} 
[\ldots, A, 0, B, \ldots] \rightarrow [\ldots, A+B, \ldots]
\eeq 
(see Proposition 3 in \cite{folded}, for instance).   
\end{rem}

In the same spirit as \cite{Hone3}, 
it is perhaps more natural to replace the first term in (\ref{eq:Pierce}) with 
an arbitrary $r\in\Q$,  resulting in a rational number minus a Pierce series. 
Such a modified Pierce series can be obtained immediately from 
(\ref{alpie1}), simply by regrouping the terms as 
\beq\label{regroup}
\left(\frac{p}{q}  +\frac{1}{x_2}\right)-\frac{1}{x_2}+
\frac{1}{x_3}-\ldots = \frac{p'}{q'}-\sum_{j=2}^\infty \frac{(-1)^j}{x_j'}, 
\eeq
where, since $q=x_1$ divides $x_2$,  $p'=p (x_2/q) +1\in\Z$, $q'=x_2$ and $x_j'=x_{j+1}$ 
for $j\geq 1$. In the proof of the preceding theorem, the initial term of the regrouped 
series appears with the continued fraction 
expansion   
$p'/q'=S_2=[a_0;a_1,\ldots, a_{2k+3}]$, ending in a partial quotient with an odd index. 

In order to formulate the most general result possible, we once again 
start with an arbitrary $r=p/q\in\Q$ and $m\in\Zo$, but this time 
take the continued fraction expansion 
\beq\label{oddcf} 
\frac{p}{q}=  [a_0;a_1,a_2,\dots,a_{2k+1}], 
\eeq  
ending in an odd index, which is always possible by (\ref{props}). 
The following result includes the regrouped series (\ref{regroup}) 
obtained from Theorem \ref{pierce1} as the special case $m=0$. 
 
\begin{thm}\label{pierce2} 
Suppose $x_1=q>1$, and let $\al$ be given by the series  
\beq \label{alpie2} 
\al = \frac{p}{q}-\sum_{j=2}^\infty \frac{(-1)^j}{x_j}, 
\eeq 
for $(x_n)$ defined by (\ref{recrels}) with 
\beq\label{y0valpie2}
y_0=mq+q_{2k}-1, 
\eeq 
where $q_{2k}$ is the denominator of the $2k$th convergent of 
(\ref{oddcf}), and the restriction $m\geq 1$ should be imposed 
if $q_{2k}=1$. 
Then  the continued fraction expansion of $\al$ has the form (\ref{al}), 
where the partial quotients $a_j$ coincide with those of 
(\ref{oddcf}) for $0\leq j\leq 2k+1$, while  
\beq\label{pie2cf} 
a_{2k+2}=m+u_1 y_0-1, 
\eeq 
and 
\beq \label{pie2coeff}
a_{2k+3j}=1, \, \, a_{2k+3j+1}=x_j-1, \, \, a_{2k+3j+2}=u_{j+1}y_j-1 \,\,\, %for \quad 
\forall j\geq 1, 
\eeq 
with the denominators of the convergents being given by 
\beq\label{cvgtspie2} 
q_{2k+3j-2}=x_j, \,\, 
q_{2k+3j-1}=y_j-x_j+1, \,\,
q_{2k+3j}=y_j+1 \,\, \, \forall j\geq 1. 
\eeq 
\end{thm} 

We omit the proof of the above,  based on showing that the partial sums are given by 
$$ 
\frac{p}{q}-\sum_{j=2}^N \frac{1}{x_j}= \frac{p_{2k+3N-2}}{q_{2k+3N-2}}, 
$$ 
since the steps 
are essentially the same 
as for Theorem \ref{pierce1}, but  the formulae (\ref{cvgtspie2})
%for the denominators of the convergents 
have been  included for completeness.

\begin{rem} Similarly to (\ref{regroup}), the case $m=0$ of Theorem \ref{pierce1} can be 
obtained from Theorem \ref{pierce2} by combining the first two terms of the series 
(\ref{alpie2}) into one. The series 
(\ref{aleng}) can be reduced to the case $m=0$ of Theorem \ref{engelthm} in 
the same way.
\end{rem}

%----------------
\section{Irrationality exponents }
\label{sec:transc}
%----------------
\setcounter{equation}{0}

The irrationality exponent $\mu(\al)$ of a real number $\al$ is defined to be the 
supremum of the set of real numbers $\mu$ such that there are  infinitely many 
rational numbers $P/Q$ satisfying the  inequality 
$0<|\al - P/Q|<1/Q^\mu$. For an irrational number, $\mu(\al)\geq 2$, since 
the convergents of its continued fraction expansion (\ref{al}) provide 
infinitely many $P/Q$ with $|\al - P/Q|<1/Q^2$. In fact, in the sense 
of Lebesgue measure, almost all real numbers 
have irrationality exponent equal to 2, while a famous theorem 
of Roth \cite{Roth} says that every algebraic irrational number 
has $\mu(\al)=2$. Large classes of transcendental numbers with $\mu(\al)=2$ 
are presented in \cite{bugeaudetal}, but as we shall see, the 
transcendental numbers defined by the modified 
Engel and Pierce series above do not belong to these classes.  

The transcendence of each of the numbers $\al$ defined  in the previous section 
is essentially a consequence of the rapid growth of the associated sequence 
$(x_n)$.  
According to the result of Lemma 2.2 in \cite{Hone3}, 
$x_{n+1}>x_n^{5/2}$ for $n\geq3$,  which allows one to show that 
$\mu(\al)\geq 5/2$, but here  
we present a significant improvement on this result. 

\begin{lem} For all $\eps>0$, there is some $N$ such that   
\beq\label{xbd} 
x_{n+1}>x_n^{\mu^* -\eps}, \qquad \mu^*=\frac{3+\sqrt{5}}{2}, 
\eeq 
for all $n\geq N$. 
\end{lem} 
\begin{prf}
From 
(\ref{nonl}) and the fact that $y_n = x_{n+1}/x_n>1$ for $n\geq 1$, 
it is clear that $x_{n+1}>x_n^2$ for $n\geq 2$, which implies that 
$x_{n-1}<x_n^{1/2}$ for $n\geq 3$. Thus,  using  (\ref{nonl}) once again 
it follows 
that $x_{n+1}>x_n^3/x_{n-1}>x_n^{5/2}$, which 
is the basic estimate given in \cite{Hone1, Hone3}. To improve on 
this, we proceed by induction, assuming that 
$$ 
x_{n+1}>x_n^{\rho_k}\quad \mathrm{for} \,\, n\geq k+2, 
$$
and then from (\ref{nonl}), the same argument as before (for $k=0$) 
gives  
\beq \label{rho} 
x_{n+1}>x_n^3/x_{n-1}>x_n^{\rho_{k+1}} 
\,\,\,\mathrm{for} \,\, n\geq k+3, \quad 
\mathrm{where} \,\,\, \rho_{k+1}=3-\rho_k^{-1},
\eeq 
with $\rho_0=2$. 
The solution of the recurrence for 
$\rho_k$ in  (\ref{rho}) is obtained via 
$\rho_k=f_{k+1}/f_k$, which implies $f_{k+2}-3f_{k+1}+f_k=0$,  
so that the sequence $(f_k)$ is just $1,2,5,13,34,\ldots$, i.e.\ 
a bisection of the Fibonacci numbers. Hence 
$\lim_{k\to\infty}\rho_k = (3+\sqrt{5})/2=\mu^*$, and the 
result follows. 
\end{prf}  

\begin{cor} All of the numbers $\al$ defined in Theorems \ref{engelthm}, 
\ref{pierce1} and \ref{pierce2} are transcendental, with irrationality 
exponent $\mu(\al)\geq (3+\sqrt{5})/2$. 
\end{cor} 
\begin{prf} This follows from the same argument  as used to prove Theorem 4
in \cite{Hone1}, and  Theorem 2.3 in \cite{Hone3}, so we only sketch the details. 
For any $\eps>0$, the terms of the sequence $(x_n)$ satisfy the inequality 
(\ref{xbd}) for sufficiently large $n$. 
Due to the fact that the partial sums 
of the series (\ref{aleng}), (\ref{alpie1}) and (\ref{alpie2}) coincide with 
particular convergents of the continued fraction expansion of $\al$, 
and the fact that, from (\ref{cvgtseng}),  (\ref{cvgtspie1}),  and  (\ref{cvgtspie2}), 
the denominators of these convergents %happen to 
coincide with the terms 
of the sequence $(x_n)$,  a comparison with a geometric series 
shows that for all $\delta>0$, the inequality 
$$ 
\left|\al -\frac{P}{Q}\right| < \frac{1}{Q^{\mu^*-\delta}}
$$  
holds for infinitely many rational approximations $P/Q$.
\end{prf} 

The irrationality exponent can be computed explicitly in terms of 
the continued fraction expansion of $\al$, using one of the formulae below: 
\beq\label{limsup} 
\mu(\al) = 1+ \limsup\limits_{n\rightarrow\infty} \frac{\log q_{n+1}}{\log q_n} 
=2+\limsup\limits_{n\rightarrow\infty} \frac{\log a_{n+1}}{\log q_n} 
\eeq 
(see \cite{bugeaudetal}, or Theorem 1 in \cite{sondow}, for instance).
If the function $\widehat{G}$ in (\ref{naut}) is chosen suitably, 
then these limits can be evaluated precisely.

\begin{thm} For some integer $d\geq 1$, let 
\beq\label{las}\la=\frac{d+2+\sqrt{d(d+4)}}{2},\eeq 
and let $\al$ be given by a modified Engel or Pierce series,   defined according to one 
of  Theorems \ref{engelthm}, 
\ref{pierce1} or \ref{pierce2}, with the sequence $(u_n)$ in (\ref{recrels}) 
being specified by a polynomial in $x_n$, namely 
\beq \label{un} 
u_n = \sum_{j=0}^{d-1} v_n^{(j)} \, x_n^{d-1-j},
\eeq 
where for each $j$,  the coefficient $(v_n^{(j)})$ is an integer sequence, with 
$ v_n^{(0)}\in \Zp$, $n\geq 1$, and as $n\to\infty$, for some $\nu<\la$, 
\beq\label{vbds} 
\log v_n^{(0)}=O(\nu^n),  \quad and \quad v_n^{(j)}=O( v_n^{(0)}), \,\, j=1,\ldots, d-1.
\eeq  
Then $\al$ has irrationality exponent $\mu(\al)=\la$. 
%be the largest root of the quadratic  
\end{thm}
\begin{prf}
If $u_n$ is given by (\ref{un}), then setting 
$\La_n=\log x_n$ and taking logarithms in (\ref{nonl}) yields
\beq\label{qualin}
\La_{n+1}-(d+2)\,\La_n +\La_{n-1}=\De_n, 
\eeq 
where, with 
$v_n^{(d)}=1$, 
$$ 
\Delta_n  = \log v_n^{(0)} + 
\log \left(1+ \sum_{j=1}^{d} \frac{v_n^{(j)}}{ v_n^{(0)}} \, x_n^{-j}\right) 
= \log v_n^{(0)} + O(x_n^{-1})=O(\nu^n)$$
as $n\to \infty$, by (\ref{vbds}).
By adapting the method of \cite{ahosloane},  the solution of (\ref{qualin}) 
is found formally as 
\beq \label{form}
\La_n = A\, \la^n + B\, \la^{-n} + \sum_{j=1}^{n-1} 
\left(\frac{\la^{n-j}-\la^{j-n}}{\la-\la^{-1}}\right) \, \De_j, 
\eeq 
for constants $A,B$ which can be fixed from the 
initial values $\La_1=\log q$, $\La_2=\log x_2$ 
(cf.\ Proposition 5 in \cite{Hone1}).  Hence 
\beq\label{asy} 
\La_n\sim C\,\la^n, \qquad C= A + \sum_{j=1}^\infty \frac{\la^{-j}\De_j}{\la-\la^{-1}}, 
\eeq 
and $C>0$ since $x_n\to\infty$ as $n\to\infty$. Then in the case 
of modified Engel series, from the formulae (\ref{cvgtseng}) there are two different 
cases for the first limit in (\ref{limsup}): when $n$ is even,  the limit 
is 
$$ 
1+\lim_{N\to\infty}\frac{\log (y_N-1)}{\log x_N}=
1+ \lim_{N\to\infty}\frac{\La_{N+1}-\La_N}{\La_N}=\la 
$$ 
from the asymptotic behaviour (\ref{asy}), while for odd $n$ it is % becomes 
$
1+1/(\la-1)\leq \la$, 
where this inequality holds for all $\la\geq(3+\sqrt{5})/2$; 
so the limit superior is $\la$. For the modified Pierce 
series in Theorem \ref{pierce1}, from  (\ref{cvgtspie1}) 
the corresponding limit varies with $n\bmod 3$, and 
gives $\la$, $2$ or $1+\la/(\la -1)\leq\la$, so again the limit superior is $\la$, 
and the case of Theorem \ref{pierce2} is identical. 
\end{prf}

\begin{thm} For all $\nu\geq(3+\sqrt{5})/2$ there are 
infinitely many $\al$ defined by modified Engel or Pierce series 
with irrationality exponent $\mu(\al)=\nu$.   
\end{thm} 
\begin{prf} If $\nu$ is one of the special values (\ref{las})  
then for any fixed $r=p/q$ and $m$ there are 
uncountably many choices of the sequence $(u_n)$ that 
take the form (\ref{un}), satisfy (\ref{vbds}) and produce $\al$ with 
distinct continued fraction expansions. So suppose that 
$d\geq 1$ (and hence $\la$) is fixed, and take any $\nu>\la$ and 
$(u_n)$ of the same form as in (\ref{un}), except that now 
\beq\label{vasy}
\log v_n^{(0)} \sim C' \, \nu^n,  \quad \mathrm{and} \,\,\, v_n^{(j)}=O( v_n^{(0)}), \,\, j=1,\ldots, d-1
\eeq  
for some $C'>0$ (for instance, one could take 
$ v_n^{(0)}=\left \lceil{\exp (C'\nu^n)}\right \rceil$). 
Then the formal expression (\ref{form}) for the solution of (\ref{qualin}) 
is still valid, but now from (\ref{vasy}) 
the sum on the right is the dominant term, growing like 
$\nu^n$ as $n\to\infty$. Substituting this leading order asymptotic behaviour 
back into (\ref{qualin}) yields 
$$ 
\La_n\sim C\, \nu^n, \qquad C=\Big(\nu-(d+2)+\nu^{-1}\Big)^{-1}\, C'>0, 
$$ 
and then the first limit 
in (\ref{limsup}) is evaluated as before to yield $\mu(\al)=\nu$. 
\end{prf}

\section{More explicit continued fractions} 
\setcounter{equation}{0}

In \cite{Hone2} one of us obtained the explicit continued fraction expansion 
for an Engel series (\ref{eq:Engel}) with %$x_1=1$ and 
the stronger 
divisibility property 
\beq\label{strong} 
x_j^2|x_{j+1},  \qquad 
\mathrm{with} \,\,\,
%positive  integers \beq\label{zseq}
z_{j+1}=\frac{x_{j+1}}{x_j^2} \in \Zp, \quad  j\geq 1
\eeq
being an arbitrary sequence of ratios. (In fact the condition $z_j\geq 2$ was imposed in \cite{Hone2},
but it was explained how to deal with some $z_j=1$ by applying the 
operation (\ref{concat}).) 
Series of this particular form include  %the sum 
\beq\label{kempner} 
\sum_{n=0}^\infty \frac{1}{u^{2^n}} 
\eeq 
for integer $u\geq 2$, which for $u=2$ is known as the Kempner number. All of the numbers (\ref{kempner}) 
are transcendental, with irrationality exponent 2 \cite{bugeaudetal}; 
their 
continued fraction expansions were  found in 
recursive form in \cite{shallit1}, with a non-recursive representation 
described in \cite{shallit2}, and further generalizations with a 
similar recursive structure being  given in \cite{shallit3} and later \cite{folded}. 
The continued fraction expansion of an Engel 
series  (\ref{eq:Engel}) %having this 
with the stronger 
divisibility property above has the same sort of 
recursive structure, defined by a particular subsequence of the 
convergents with finite continued fractions  whose length approximately 
doubles at each step. Here we 
further generalize this result by considering an arbitrary 
$r=p/q\in\Q$ added to an Engel series of this type, 
of the form (\ref{pqengel}) with the property  (\ref{strong}). 

%In order to present the more general result
Given a finite continued fraction (\ref{ficf}) for $n\geq 2$, written  as 
$
[a_0; {\bf a}]$, where ${\bf a}=(a_1,a_2,\ldots,a_n)$ is the word 
defining the fractional part, it is 
convenient to define the following family of transformations:
\beq\label{zfam} 
\varphi_z : \quad [a_0; {\bf a}] \mapsto [a_0; {\bf a}, z-1, \hat{{\bf a}}], 
 \eeq where the word $\hat{{\bf a}}$ is given by $$ 
%\quad \mathrm{where} \,\,\,
\hat{{\bf a}}=(1,a_n-1,a_{n-1},\ldots,a_2,a_1). 
$$ 
For each $z$, starting from   a continued fraction whose 
final partial quotient has index $n$, 
$\varphi_z$ produces a new continued fraction
$\varphi_z([a_0;{\bf a}])$ whose 
final partial quotient has index $2n+2$. It is also 
helpful to define the one-parameter family of transformations 
\beq\label{zsfam} 
\varphi_z^* : \quad [a_0; {\bf a}] \mapsto [a_0; {\bf a}, z-1, {\bf a}^*], 
 \eeq where the word at the end is  
$$ 
%\quad \mathrm{where} \,\,\,
{\bf a}^*=(a_{n-1}+1,a_{n-2},\ldots,a_2,a_1). 
$$ 
For the latter family, $[a_0;{\bf a}]$ is sent to 
 $\varphi_z^*([a_0;{\bf a}])$, whose 
final partial quotient has index $2n$. 
The transformation $\varphi^*_z$ is just the result of applying 
$\varphi_z$ to a continued fraction (\ref{ficf}) for which $a_n=1$, 
followed by using the concatenation operation (\ref{concat}) 
to remove the zero that appears. 

Our interest in the above transformations is due to 

\begin{lem}\label{trans} 
$$ \frac{p_n}{q_n}+ \frac{(-1)^n}{zq_n^2} = \varphi_z([a_0;{\bf a}]). $$  
\end{lem}
\begin{prf} This is a corollary of Proposition 2 in \cite{folded}, and 
follows from the remarks made beneath the proof given  there. 
The case of even $n=2k$
%, which is the only case of interest in the sequel, 
also follows from the proof of 
Proposition 2.1 in \cite{Hone2}. 
\end{prf}

We now consider a rational number $r=p/q$ added to an Engel series 
with the property (\ref{strong}).  Without loss of generality, we 
exclude the case $r\in \Z$ (when $q=x_1=1$), since when 
$x_2\geq 2$ one can instead take $p'/q'=p/q+1/x_2$ 
as the initial rational term. 

\begin{thm}\label{modstreng}
Let $p/q\in\Q\setminus\Z$ be a rational number in lowest 
terms, with continued fraction expansion 
$$ 
\mathcal{C}_1 = [a_0; {\bf a}] 
$$ 
taken in the form (\ref{evencf}), and let $(x_n)_{n\geq 1}$ be the integer 
sequence defined by 
$$x_1=q, \qquad x_{j+1}=x_j^2 z_{j+1}, \quad j\geq 1, $$ 
where $(z_n)_{n\geq 2} $ is an arbitrary 
sequence of positive integers. Define a 
sequence of finite 
continued fractions according to 
$$ 
\mathcal{C}_2 = 
\begin{cases}
  \varphi_{z_2}^*(\mathcal{C}_1), & \text{if $a_{2k}=1$},\\
   \varphi_{z_2}(\mathcal{C}_1), & \text{otherwise},
\end{cases}
$$  
and, for all $j\geq 2$, with  
 the restriction $z_j\geq 2$ imposed, 
$$ 
\mathcal{C}_{j+1} = 
\begin{cases}
  \varphi_{z_{j+1}}^*(\mathcal{C}_j), & \text{if $a_{1}=1$},\\
   \varphi_{z_{j+1}}(\mathcal{C}_j), & \text{otherwise}.
\end{cases}
$$ 
Then the modified Engel series 
\beq \label{alstreng} 
\al = \frac{p}{q}+\sum_{j=2}^\infty \frac{1}{x_j}
\eeq 
has continued fraction expansion $\lim_{N\to\infty}\mathcal{C}_N$. 
\end{thm} 

\begin{prf} We show by induction 
that the partial sum 
$$ 
S_N = \frac{p}{q}+\sum_{j=2}^N \frac{1}{x_j}
$$ has continued fraction expansion $\mathcal{C}_N$, 
and then the result for $\al$ given by (\ref{alstreng}) 
follows in the limit $N\to\infty$. 
When $N=1$ we have 
$
S_1=%\frac{p}{q}=\frac{p_{n_1}}{q_{n_1}}
p/q=p_{n_1}/q_{n_1}=\mathcal{C}_1
$,
where $n_1=2k$ is the index of the final 
partial quotient of $\mathcal{C}_1$, chosen to be even 
as in (\ref{evencf}). 
If $a_{n_1}=a_{2k}=1$, then we apply 
the transformation (\ref{zsfam}) with $z=z_2$, 
and otherwise we apply (\ref{zfam}). 
In either case we obtain a new continued fraction 
$\mathcal{C}_2$ whose final partial quotient has an even index, 
$n_2$ say, with $a_{n_2}=a_1$. 
 At each subsequent step $N=j$, 
we have a continued fraction $\mathcal{C}_j$, 
with an even number $n_j$ being the index 
of its final partial quotient, that is $a_{n_j}=a_1$, 
so that to obtain  $\mathcal{C}_{j+1}$,  if $a_1=1$ then
we must apply 
(\ref{zsfam}) with $z=z_{j+1}$, 
or  (\ref{zfam}) otherwise. 
Then, from Lemma \ref{trans} 
we have
$$ 
S_{j+1}=S_j+\frac{1}{x_{j+1}}
= \frac{p_{n_j}}{q_{n_j}}+\frac{1}{z_{j+1}q_{n_j}^2}
=\mathcal{C}_{j+1},
$$
where  we used the additional hypothesis that $q_{n_j}=x_j$.  
Then the 
induction is almost complete, apart from the verification that 
$q_{n_{j+1}}=x_{j+1}$ is also a consequence of the above. However,  
$\mathcal{C}_{j+1}$ equals 
$$ 
 \frac{p_{n_{j+1}}}{q_{n_{j+1}}} = \frac{ p_{n_j} (x_{j+1}/q_{n_j})+1}
{x_{j+1}} =\frac{p_{n_j}x_jz_{j+1}+1}{x_{j+1}}, 
$$ 
since $q_{n_j}=x_j$ and  $x_{j+1}=x_j^2z_{j+1}$. Then any prime $P$ that divides $x_{j+1}$ must divide 
$x_j$ or $z_{j+1}$, so  the numerator on the right-hand side above is 
congruent to $1\bmod P$. Thus the fraction on the right-hand side, with denominator $x_{j+1}$, is 
in lowest terms, and since the convergent  $p_{n_{j+1}}/q_{n_{j+1}}$ 
is also in lowest terms this means that $q_{n_{j+1}}=x_{j+1}$ as required. 
\end{prf} 

\begin{rem}
The case where one or more of the $z_j=1$  (as in the sum (\ref{kempner}), 
for instance) can be dealt with 
by applying (\ref{concat}). For various examples of this, see \cite{Hone2}. 
\end{rem}

\begin{exa} For a given denominator $q$ it is sufficient to consider $0<p/q<1$, so 
 picking $q=x_1=5$ and $z_j=(j+1)^2+1$ for $j\geq 2$, $\al$ is the sum 
$$ 
%\al= 
\frac{p}{5}+ \frac{1}{5^2 \cdot 10}+\frac{1}{5^4\cdot 10^2\cdot 17}
+\frac{1}{5^8\cdot 10^4\cdot 17^2\cdot 26}
+\frac{1}{5^{16}\cdot 10^8\cdot 17^4\cdot 26^2\cdot 37}
+
\cdots ,
$$ 
and we  list the continued fractions for $1\leq p\leq 4$ in order: 
%$$ \begin{array}{l} 
%p=1 
%[0; 4,1,9,5,16,1,4,9,1,4,25,1,3,1,9,4,1,16,5,9,1,4,36,\ldots]  \\ \relax
%p=2  
%[0; 2,2,9,1,1,2,16,1,1,1,1,9,2,2,25,1,1,2,9,1,1,1,1,16,2,1,1,9,2,2,%36,
%\ldots] \end{array} $$ 
$$ 
[0; 4,1,9,5,16,1,4,9,1,4,25,1,3,1,9,4,1,16,5,9,1,4,36,1,3,1,9,5,16,\ldots]  ,
$$ 
$$  
[0; 2,2,9,1,1,2,16,1,1,1,1,9,2,2,25,1,1,2,9,1,1,1,1,16,2,1,1,9,2,2,%36,
\ldots] ,
$$
$$ 
[0;1,1,1,1,9,2,1,1,16,2,2,9,1,1,1,1,25,2,1,1,9,2,2,16,1,1,2,9,1,1,%1,
\ldots ],
$$
$$ 
[0;1,4,9,1,3,1,16,4,1,9,4,1,25,5,9,1,4,16,1,3,1,9,4,1,36,5,9,1,3,%1,
\ldots].
$$ 
\end{exa} 

\small
\noindent
{\bf Acknowledgments:} ANWH is funded by EPSRC Fellowship EP/M004333/1, 
and thanks the School of Mathematics \& Statistics, UNSW for 
hospitality and support under the Distinguished Researcher Visitor Scheme.
%from the 
%Engineering and Physical Sciences Research Council. 
JLV is supported by Grant MTM2015-65888-C4-4-P from MINECO/FEDER.

%----------------

%----------------

%----------------
\end{document}